\documentclass[12pt]{article}
\usepackage[utf8]{inputenc}
\usepackage{amsmath, amssymb, amsthm, enumerate} 

\newtheorem{theorem}{Theorem}
\newtheorem{proposition}[theorem]{Proposition}
\newtheorem{definition}[theorem]{Definition}
\newtheorem{lemma}[theorem]{Lemma}
\newtheorem{corollary}[theorem]{Corollary}

\newtheorem{remarks}[theorem]{Remarks}
\newtheorem{remark}[theorem]{Remark}

\newtheorem{notation}[theorem]{Notation}

\def\am{^{-1}}

\def\cc{{\mathbb C}}

\def\su{\subset}
\def\sp{\supset}

\def\al{\alpha}
\def\be{\beta}

\def\Ga{\Gamma}

\def\ep{\varepsilon}

\def\ph{\phi}

\def\la{\lambda}
\def\La{\Lambda}

\def\cd{\cdot}
\def\stb{,\ldots ,}

\def\del{\partial}

\def\V{\Vert}

\def\({\left(}
\def\){\right)}
\def\[{\left[}
\def\]{\right]}

\def\ol{\overline}
\def\ds{\displaystyle}

\def\sumir{\sum_{i=1}^r}

\def\sumi0n{\sum_{i=0}^n}

\def\x1n{x_1 \stb x_n}
\def\y1n{y_1 \stb y_n}
\def\deg{{\rm deg}\, }
\def\dim{{\rm dim}\, }

\def\cl{{\rm cl}\, }

\def\noi{\noindent}

\begin{document}

\title{Translation invariant linear spaces of polynomials}

\author{Gergely Kiss and Mikl\'os Laczkovich (Budapest, Hungary)}

\footnotetext[2]{{\bf MR subject classification:} 32A08, 13C05}
\footnotetext[3]{Both authors were supported by the Hungarian National
Foundation for Scientific Research, Grant No. K124749. The first author
was also supported by the Premium Postdoctoral Fellowship of the Hungarian
Academy of Sciences.}

\maketitle

\begin{abstract} 
A set of polynomials $M$ is called a {\it submodule} of $\cc [x_1 \stb x_n ]$
if $M$ is a translation invariant linear subspace of $\cc [x_1 \stb x_n ]$. We 
present a description of the submodules of $\cc [x,y]$ in terms of a
special type of submodules. We say that the submodule $M$ of $\cc [x,y]$ is an
{\it L-module of order} $s$ if, whenever $F(x,y)=\sum_{n=0}^N f_n (x)
\cd y^n \in M$ is such that $f_0 =\ldots = f_{s-1}=0$, then $F=0$. 
We show that the proper submodules of $\cc [x,y]$ are the
sums $M_d +M$, where $M_d =\{ F\in \cc [x,y] \colon \deg _x F <d\}$, and $M$ is
an L-module. We give a construction of L-modules parametrized by sequences of complex numbers.

A submodule $M\su \cc [x_1 \stb x_n ]$ is {\it decomposable} if it is
the sum of finitely many proper submodules of $M$. Otherwise $M$ is
{\it indecomposable}. It is easy to see that every submodule of $\cc [x_1 \stb
x_n ]$ is the sum of finitely many indecomposable submodules. In $\cc [x,y]$
every indecomposable submodule is either an L-module or equals $M_d$
for some $d$. In the other direction we show that 
$M_d$ is indecomposable for every $d$, and so is every L-module
of order $1$.

Finally, we prove that there exists a submodule of $\cc [x,y]$ (in fact, an
L-module of order $1$) which is not relatively closed in $\cc [x,y]$.
This answers a problem posed by L. Sz\'ekelyhidi in 2011.
\end{abstract}

\section{Introduction and main results}

In this note we are concerned with the translation invariant linear subspaces of
$\cc [x_1 \stb x_n ]$, the ring of polynomials of $n$ variables having complex
coefficients. By making use of Taylor's formula it is not difficult to see that
a linear subspace of $\cc [x_1 \stb x_n ]$ is translation invariant if and only
if it is invariant under partial differentiation \cite[Lemma 7]{LSz}. Thus a
linear subspace of $\cc [x_1 \stb x_n ]$ is translation invariant if and only
if it is a module under the ring of partial differential operators. For this
reason we say that $M$ is a {\it submodule of} $\cc [x_1 \stb x_n ]$ (or briefly
a {\it module}) if $M$ is a translation invariant linear subspace of
$\cc [x_1 \stb x_n ]$. 

It is easy to check that the only submodules of $\cc [x]$ are $\cc [x]$
itself and the modules $\{ f\in \cc [x] \colon \deg f<d\}$ $(d=0,1,\ldots )$.

Simple examples of submodules of $\cc [x,y]$ are $\cc [x,y]$ itself, $\cc [x]$,
$\cc [y]$, $\{ f(x+y)\colon f\in \cc [x] \}$, $\{ f(ax+by)\colon f\in \cc [x]
\}$ $(a,b\in \cc )$, $\{ f(x)+g(y)\colon f,g\in \cc [x] \}$, $\{ f\in \cc [x,y]
\colon \deg f <d\}$ $(d=0,1,\ldots )$, $\{ f\in \cc [x,y]\colon \deg_x f <d_1 ,
\ \deg_y f <d_2 \}$ $(d_1 ,d_2 =0,1, \ldots )$. Here $\deg _x f$ and $\deg _y f$ 
denote the degree in the variable $x$ (resp. $y$) of the polynomial $f\in
\cc [x,y]$.

Each of these modules is relatively closed in $\cc [x,y]$ in the following
sense: if $f_n$ belongs to the module $M$ in question for every $n$ and
$f_n \to f\in \cc [x,y]$ pointwise (or uniformly on compact sets), then
$f\in M$. 

The investigations of these note were motivated by the following problem posed
by L. Sz\'ekelyhidi \cite{Sz}: {\it is it true that every submodule of
$\cc [x_1 \stb x_n ]$ is relatively closed in $\cc [x_1 \stb x_n ]$?}
In other words, is every submodule of $\cc [x_1 \stb x_n ]$ a variety?
In Theorem \ref{t6} we show that the answer to Sz\'ekelyhidi's question is
negative. Our example is a special case of a general construction of some
submodules of $\cc [x,y]$, called L-modules.

We represent the elements of $\cc [x,y]$ in the form
\begin{equation}\label{e1}
F(x,y)=\sum_{n=0}^\infty f_n (x)\frac{y^n}{n!},
\end{equation}  
where $f_n \in \cc [x]$ for every $n$, and $f_n =0$ if $n$ is large enough.
We say that the module $M\su \cc [x,y]$ is an {\it L-module of order $s$} if,
whenever $F$ in \eqref{e1} belongs to $M$ and such that $f_n =0$ for every
$n<s$, then $F=0$. In Section 2 we give a construction of L-modules parametrized
by sequences of complex numbers (Theorem \ref{t1}).

Let $M_d$ denote the module $\{ f\in \cc [x,y] \colon \deg _x f <
d\}$. In Section 3 we show that every proper submodule of $\cc [x,y]$
can be represented in the form $M_d +M$, where $M$ is an L-module (Theorem 
\ref{t2}). Under some mild restrictions on $M$, the representation is unique
(see Remark \ref{r1}). The obstacles in the way of generalizing this result for
polynomials of more than two variables are discussed in Remark \ref{r2}.

We say that a submodule of $\cc [x_1 \stb x_n ]$ is indecomposable, if it
cannot be written as a finite sum of proper submodules.
It is easy to see that every submodule of $\cc [x_1 \stb x_n ]$ is the sum
of finitely many indecomposable submodules (Proposition \ref{p1}). It follows
from Theorem \ref{t2} that every indecomposable submodule of $\cc [x,y]$ is
either an L-module or equals $M_d$ for some $d$. In the other direction
we prove that $M_d$ is indecomposable for every $d$, and that all L-modules
of order $1$ are indecomposable submodules of $\cc [x,y]$ (Theorems \ref{t4} and
\ref{t5}).

\section{L-modules}
Let $S=\{ f(x+y)\colon f\in \cc [x] \}$. It is clear that $S$ is a submodule
of $\cc [x,y]$. Since $f(x+y)=\sum_{n=0}^\infty f^{(n)} (x)\cd \tfrac{y^n}{n!}$
by Taylor's formula, it follows that the elements of $S$ are the polynomials
$\sum_{n=0}^\infty f_n \cd \tfrac{y^n}{n!}$, where $f_n \in \cc [x]$ for every $n$,
and $f_n =f'_{n-1}$ for every $n\ge 1$. In particular, $S$ has the property that
if $F=\sum_{n=0}^\infty f_n \cd \tfrac{y^n}{n!} \in S$ and $f_1 =0$, then $F=0$.

The module $S$ is the prototype of the submodules we are about to define.
\begin{notation}
{\rm Every polynomial $F\in \cc [x,y]$ can be represented uniquely in the
form \eqref{e1}, where $f_0 ,f_1 ,\ldots \in \cc [x]$, and $f_n =0$ if $n$
is large enough. The polynomials $f_n$ will be called} the coordinate
polynomials {\rm of $F$, and will be denoted by $[F]_n$ $(n=0,1,\ldots )$.

If $A\su \cc [x,y]$ and $s$ is a positive integer, then we put}
$$V_{A,s} =\{ ([F]_0 \stb [F]_{s-1}) \colon F\in A\} .$$
\end{notation}

Clearly, if $A$ is a module, then $V_{A,s}$ is a linear subspace of $\cc [x]^s$.
Note that if $F\in A$, then $([F]_{k-s}\stb [F]_{k-1})\in V_{A,s}$ for every
$k\ge s$. This follows from $\tfrac{\del ^{k-s}}{\del y^{k-s}} F\in A$.
The set $V_{A,s}$ also has the following property: if $(f_0 \stb f_{s-1})\in
V_{A,s}$, then $(f'_0 \stb f'_{s-1})\in V_{A,s}$. This is clear from the fact that
if the polynomial in \eqref{e1} belongs to $A$, then
\begin{equation}\label{e16}
\frac{\del}{\del x} F=\sum_{n=0}^\infty f'_n (x)\frac{y^n}{n!} \in A.
\end{equation}  

\begin{definition} \label{d1}
{\rm Let $s$ be a positive integer. We say that $M\su \cc [x,y]$ is an}
$L$-module of order $s$ {\rm if $M$ is a module and, whenever $F\in M$ and
$[F]_n =0$ for every $n<s$, then $F=0$.}
\end{definition}
Since the submodules of $\cc [x,y]$ are linear spaces, the condition formulated
in Definition \ref{d1} is equivalent to the following: if $F\in M$, then
$F$ is determined by the coordinate polynomials $[F]_0 \stb [F]_{s-1}$.

It is clear from the definition that if $M$ is an L-module of order $s$, then
it is also an L-module of order $t$ for every $t\ge s$.

\begin{theorem} \label{t1}
Let $M\su \cc [x,y]$ be an $L$-module of order $s$. Then there exists a linear
map $L\colon \cc [x]^s \to \cc [x]$ such that for every $F\in M$ we have
\begin{equation}\label{e2}
[F]_n = L([F]_{n-s} \stb [F]_{n-1} )
\end{equation}  
for every $n\ge s$. 
More precisely, there are complex numbers $a_{i,j}$ $(i=1\stb s, \
j=1,2,\ldots )$ such that \eqref{e2} holds for every $n\ge s$, where 
\begin{equation}\label{e3}
L (f_1 \stb f_s )=\sum_{i=1}^s \sum_{j=1}^\infty a_{i,j} f_i^{(j)}
\end{equation}  
for every $f_1 \stb f_s \in \cc [x]$. (Note that the sum in the right hand side
of \eqref{e3} only has a finite number of nonzero terms for every $f_1 \stb f_s
\in \cc [x]$).
\end{theorem}

\noi
{\it Proof.} If $F\in M$, then we put $L([F]_0 \stb
[F]_{s-1})=[F]_s$. This definition makes sense, since $M$ is an L-module, and
thus $[F]_s$ is uniquely determined by $[F]_0 \stb [F]_{s-1}$. In this way we
defined $L$ on the set $V_{M,s}$. It is clear that $L$ is linear. 

Suppose $F\in M$, and let $k\ge s$ be given. Then
\begin{equation}\label{e15}
\frac{\del ^{k-s}}{\del y^{k-s}}F (x,y)=\sum_{n=0}^\infty [F]_{n+k-s} (x)
\frac{y^n}{n!}   \in M,
\end{equation}
and thus $[F]_k =L([F]_{k-s}\stb [F]_{k-1})$. This proves the first statement of
the theorem including \eqref{e2}, except that $L$ is only defined on $V_{M,s}$.

If $F\in M$, then \eqref{e16} holds, and thus
$$L([F]_0 \stb [F]_{s-1} )' =L([F]'_0 \stb [F]'_{s-1} ) ,$$
since both sides equal $[F]'_s$. Therefore,
\begin{equation}\label{e4}
L(f_0 \stb f_{s-1} )' =L(f'_0 \stb f'_{s-1} )
\end{equation}
holds for every $(f_0 \stb f_{s-1} )\in V_{M,s}$. Next we prove that
\begin{equation}\label{e5}
\deg L(f_0 \stb f_{s-1} ) \le \max_{0\le i\le s-1} \deg f_i
\end{equation}
for every $(f_0 \stb f_{s-1} )\in V_{M,s}$. Indeed, let $F\in M$,
and let $\ds{\max_{0\le i\le s-1}} \deg [F]_i =d$. Then 
$$\frac{\del ^{d+1}}{\del x^{d+1}}F (x,y)=\sum_{n=0}^\infty [F]_n^{(d+1)}
\frac{y^n}{n!} \in M.$$
Now we have $[F]_i^{(d+1)} =0$ for every $i<s$, and thus $[F]_s^{(d+1)} =0$, 
since $M$ is an $L$-module of order $s$. This proves \eqref{e5}.
Then it follows from \eqref{e2} that if $F\in M$, then $\deg [F]_n \le
\max_{0\le i\le s-1} \deg [F]_i$ for every $n$. In particular, if the coordinate
polynomials $[F]_0 \stb [F]_{s-1}$ are constants, then $[F]_n$ is constant for
every $n$.

Let $W_0$ denote the set of $s$-tuples $(c_0 \stb c_{s-1})\in \cc ^s$ such that
$(c_0 \stb c_{s-1})\in V_{M,s}$. Clearly, $W_0$ is a linear subspace of
$\cc ^s$. Let $\dim W_0 =r$, and let $(c_{i,0} \stb c_{i,s-1})$ $(i=1\stb r)$ be a
basis of $W_0$. Let $F_i (y)=\sum_{n=0}^{\infty} c_{i,n} \frac{y^n}{n!} \in
M\cap \cc [y]$ for every $i=1\stb r$. Then every element of $M\cap \cc [y]$ is
the linear combination of the functions $F_i$. Indeed, if
$$F(y)=\sum_{n=0}^\infty c_n \frac{y^n}{n!} \in M\cap \cc [y] ,$$
then there is a linear combination $\ol F$ of $F_1 \stb F_{r}$ such that
$\ol F (y)=\sum_{n=0}^\infty d_n \tfrac{y^n}{n!}$, where $d_n =c_n$ for every
$n\le s-1$. Since $\ol F \in M$ and $M$ is an $L$-module of order $s$,
it follows that $d_n =c_n$ for every $n$, and $\ol F =F$.

Thus the dimension of $M\cap \cc [y]$ is at most $r$.
Since $\dim W_0 =r$, we have $\dim (M\cap \cc [y] )=r$. Now $M\cap \cc [y]$
is a proper submodule of $\cc [y]$, and thus there is a $p\ge 0$ such that
$M\cap \cc [y] =\{ f\in \cc [y]\colon \deg f <p\}$. Clearly, we must have $p=r$.

This implies $p=r\le s$, and thus $\deg F<s$ for every $F\in M\cap \cc [y]$.
That is, $c_s =0$ whenever $\sum_{n=0}^\infty c_n \tfrac{y^n}{n!} \in
M\cap \cc [y]$. Therefore, we have $L(c_0 \stb c_{s-1})=0$ for every
$(c_0 \stb c_{s-1})\in W_0$.

We construct the numbers $a_{i,j}$ with the property that, for every $d\ge 1$,
\begin{equation}\label{e6}
L(f_0 \stb f_{s-1})=\sum_{i=0}^{s-1} \sum_{j=1}^{d-1} a_{i,j} f_i^{(j)}
\end{equation}  
whenever $(f_0 \stb f_{s-1}) \in V_{M,s}$ and $\deg f_i <d$ $(i=0\stb s-1)$.
Note that \eqref{e6} is true for $d=1$. Indeed, $\deg f_i <1$ means that $f_i$
is constant, and thus the left hand side of \eqref{e6} is zero, and so is the
right hand side, since the sums $\sum_{j=1}^{d-1}$ are empty.

Let $d\ge 1$, and suppose we have defined the numbers and
$a_{i,j}$ $(i=0\stb s-1, \ j=1\stb d-1)$ such that \eqref{e6} holds for every
$(f_0 \stb f_{s-1}) \in V_{M,s}$ and $\deg f_i <d$ $(i=0\stb s-1)$.

If $(f_0 \stb f_{s-1} ) \in V_{M,s}$ and $\deg f_i \le d$ $(i=0\stb s-1)$, then,
by \eqref{e4} and \eqref{e6},
$$(L(f_0 \stb f_{s-1} ) )'=L(f'_0 \stb f'_{s-1} )=
\sum_{i=0}^{s-1} \sum_{j=1}^{d-1} a_{i,j} f_i^{(j+1)} ,$$
and thus
$$L(f_0 \stb f_{s-1} ) =
\sum_{i=0}^{s-1} \sum_{j=1}^{d-1} a_{i,j} f_i^{(j)} + C(f_0 \stb f_{s-1}) ,$$
where $C(f_0 \stb f_{s-1})$ is constant. Clearly, the map
$$(f_0 \stb f_{s-1}) \mapsto C(f_0 \stb f_{s-1})$$
is linear. Let $f_i = \sum_{\nu =0}^d \al_{i,\nu} x^\nu$ $(i\le s-1)$. Then
$C(f_0 \stb f_{s-1} )$ only
depends on the coefficients $\al _{i,d}$. Indeed, if $(g_0 \stb g_{s-1} )\in V_{M,s}$,
where $g_i =\sum_{\nu =0}^d \be_{i,\nu} x^\nu$ and $\al _{i,d}=\be _{i,d}$
$(i\le s-1)$, then $\deg (f_i -g_i )<d$, and $C(f_0 -g_0 \stb f_{s-1} -g_{s-1} )
=0$ by \eqref{e6}. Then it follows that there are numbers $b_{i,d}$
$(i=0\stb s-1)$ such that
$$C(f_0 \stb f_{s-1} )=\sum_{i=0}^{s-1} b_{i,d}\cd \al _{i,d} =
\sum_{i=0}^{s-1} b_{i,d}\cd \frac{f_i^{(d)}}{d!} .$$
Putting $a_{i,d} =b_{i,d}/d!$, we obtain \eqref{e6} with $d+1$ in place of $d$
for every $(f_0 \stb f_{s-1} ) \in V_{M,s}$, $\deg f_i \le d$ $(i=0\stb s-1)$.
In this way we obtain the numbers $a_{i,j}$ by induction on $j$. It is clear that
the numbers $a_{i,j}$ defined above satisfy \eqref{e3} for every $(f_1 \stb f_s )
\in V_{M,s}$. Now, the right hand side of \eqref{e3} makes sense for every
$(f_1 \stb f_s )\in \cc [x]^s$, and defines a linear extension of $L$ to
$\cc [x]^s$.  This completes the proof of the theorem. \hfill $\square$

\begin{remark}
{\rm
Let $M$ be an L-module of order $s$. Then $V_{M,s}$ and $L$ are connected by
the following necessary condition:
if $(f_0 \stb f_{s-1})\in V_{M,s}$ and the sequence of polynomials is defined
by $f_n =L(f_{n-s} \stb f_{n-1})$ for every $n\ge s$, then $(f_{n-s} \stb f_{n-1})
\in V_{M,s}$ for every $n\ge s$. Indeed, let $F(x,y)=\sum_{n=0}^\infty g_n (x)
\tfrac{y^n}{n!} \in M$ be such that $g_i =f_i$ for every $i<s$. Since
$g_n =L(g_{n-s} \stb g_{n-1})$ for every $n\ge s$, it follows that $g_n =f_n$
for every $n$. For every $k\ge s$ we have \eqref{e15}, hence
$(f_{k-s} \stb f_{k-1}) \in V_{M,s}$.

In the constructions of L-modules this condition should be taken into account.
Consider the following example. Let $s=2$, $V =\{ (f,f)\colon f\in \cc [x] \}$,
and let $L$ be the identically zero map from $\cc [x]^2$ into $\cc [x]$.
Then $(f,f)\in V$, $L(f,f)=0$, but $(f,0)\notin V$ if $f\ne 0$. Accordingly,
the set $M$ of functions of the form \eqref{e1} such that $(f_0 ,f_1 )\in V$
and $f_n =L(f_{n-2},f_{n-1})$ for every $n\ge 2$ is not a module. Indeed,
$F=f(x)\cd (1+y) \in M$ for every $f\in \cc [x]$, but $\tfrac{\del}{\del y} F=
f(x)\notin M$ if $f\ne 0$.
}
\end{remark}

Note that the necessary condition above is automatically satisfied if $V_{M,s}
=\cc [x]^s$. Therefore, the following construction always produces L-modules.
\begin{notation}
{\rm 
Let $\Ga =\{ a_{i,j} \colon i=1\stb s, \ j=1,2,\ldots \}$
be a set of complex numbers, and let $M_{\Ga}$ denote the set of polynomials of
the form \eqref{e1} such that $f_n =L (f_{n-s} \stb f_{n-1})$ for every $n\ge s$,
where $L$ is defined by \eqref{e3}.

The definition of $L$ implies that
\begin{equation}\label{e7}
\deg f_k < \max_{k-s \le i\le k-1} \deg f_i 
\end{equation}
for every $k\ge s$. Therefore, we have $f_n =0$ for every
$n>s+\max_{0\le i\le s-1} \deg f_i$.}
\end{notation}

\begin{lemma}\label{l1}
$M_{\Ga}$ is an L-module of order $s$.
\end{lemma}

\noi
{\it Proof.} It is enough to show that $M_{\Ga}$ is a module. Since 
$L$ is a linear map, $M_{\Ga}$ is a linear subspace of $\cc [x,y]$. If $F(x,y)$
is defined by \eqref{e1}, then
$$\frac{\del}{\del x} F(x,y)=\sum_{n=0}^\infty f'_n (x)\cd \frac{y^n}{n!}.$$
Since $f'_n =L(f'_{n-s} \stb f'_{n-1})$ for
every $n\ge s$, we have $\tfrac{\del}{\del x} F \in M_\Ga$. We also have
$$\frac{\del}{\del y} F(x,y)=\sum_{n=1}^\infty f_n (x)\cd
\frac{y^{n-1}}{(n-1)!} =\sum_{n=0}^\infty f_{n+1} (x)\cd \frac{y^n}{n!} .$$
It is clear that $\tfrac{\del}{\del y} F \in M_\Ga$, and thus $M_\Ga$ is a
module. \hfill $\square$

\section{A representation of the proper submodules of $\cc [x,y]$}

By the sum of the sets of polynomials $A,B \su \cc [x_1 \stb x_n ]$ we mean the
set $A+B=\{ f+g\colon f\in A,\ g\in B\}$. It is easy to see that if
$A,B$ are submodules of $\cc [x_1 \stb x_n ]$, then so is $A+B$.

\begin{notation}
{\rm For every nonnegative integer $d$ we denote by $M_d$ the set of polynomials
$F\in \cc [x,y]$ such that $\deg [F]_n <d$ for every $n$.}
\end{notation}
It is easy to check that $M_d$ is a submodule of $\cc [x,y]$ for every
nonnegative integer $d$. Note that $M_0 =\{ 0\}$ and $M_1 =\cc [y]$. 

In this section our aim is to prove the following.
\begin{theorem} \label{t2}
Let $A$ be a proper submodule of $\cc [x,y]$. Then there are integers $d\ge 0$
and $s\ge 1$ such that $A= M_d +M$, where $M$ is an $L$-module of order $s$.
\end{theorem}

As for the uniqueness of the representation see Remark \ref{r1}.
First we show that the sum of an L-module and $M_d$ is always a proper
submodule of $\cc [x,y]$.

\begin{lemma}\label{l2}
If $M$ is an $L$-module of order $s$ and $d\ge 0$, then $x^d y^s \notin M_d +M$.
Consequently, $M_d +M$ is a proper submodule of $\cc [x,y]$.
\end{lemma}

\noi
{\it Proof.} Suppose $x^d y^s =F+G$, where $F\in M_d$ and $G\in M$. Then we have
\eqref{e1}, where $\deg f_n <d$ for every $n$. Thus
\begin{align*}
-G(x,y) &=F(x,y) -x^d y^s\\
&  =\sum_{n=0}^{s-1} f_n (x)\cd \frac{y^n}{n!} + (f_s (x)-s! \cd x^d ) \cd
\frac{y^s}{s!} +\sum_{n=s+1}^\infty f_n (x)\cd \frac{y^n}{n!},
\end{align*}
and $-G\in M$. By \eqref{e5} we obtain
$$d=\deg ((f_s (x)-s! \cd x^d )\le \max_{0 \le i\le s-1} \deg f_i <d,$$
a contradiction. \hfill $\square$

\begin{corollary} \label{c1}
If $A ,B$ are L-modules and $M_{d_1} +A=M_{d_2} +B$, then $d_1 =d_2$.
\end{corollary}

\noi
{\it Proof.} Suppose $d_1 <d_2$. By Lemma \ref{l2}, there is an $s$ such that
$x^{d_1} y^s \notin M_{d_1} +A=M_{d_2} +B$. However, $x^{d_1} y^s \in M_{d_2} \su
M_{d_2} +B$, which is a contradiction.\hfill $\square$

The rest of the section is devoted to the proof of Theorem \ref{t2}.

Since $A \subsetneqq \cc [x,y]$ and $A$ is a linear space, we have $x^m y^s
\notin A$ for some $m,s\ge 0$. Let $d$ be the smallest nonnegative integer such
that $x^d y^s \notin A$ for some $s\ge 0$.
Then we have $M_d \su A$. Let $s$ be the smallest nonnegative integer such that
$x^d y^s \notin A$. In the course of the proof we fix the module $A$ and the
nonnegative integers $d$ and $s$ with these properties.
\begin{lemma}\label{l4}
For every polynomial $F\in A$ we have
\begin{equation}\label{e24}
\deg [F]_n <\max (d, \max_{0\le i \le s-1} \deg [F]_i )
\end{equation}
for every $n\ge s$.
\end{lemma}

\noi
{\it Proof.} Let $e=\max_{0\le i \le s-1} \deg [F]_i$, and suppose that $m\ge
\max (d,e)$, where $m=\max_{n\ge s} \deg [F]_n$. (Note that $[F]_n =0$ if $n$
is large enough.)

First we suppose $e\le d$; then $m\ge d$. Turning to the polynomial
$\tfrac{\del ^{m-d}}{\del x^{m-d}}F$ we may assume that $m=d$. Let $k$ be the
largest index with $\deg [F]_k =d$. Then $k\ge s$, $\deg [F]_n \le d$ for every
$n$, and $\deg [F]_n < d$ for every $n>k$. Turning to the polynomial
$\tfrac{\del ^{k-s}}{\del y^{k-s}}F$ we may assume that $k=s$. Then $\deg [F]_n
\le d$ for every $n$, $\deg [F]_s = d$, and $\deg [F]_n < d$ for every $n>s$.
Since $M_d \su A$ and $x^d y^n  \in A$ for every $n<s$ by the choice of $s$, it
follows that $[F]_n \cd \tfrac{y^n}{n!} \in A$ for every $n\ne s$, and thus
$[F]_s \cd \tfrac{y^s}{s!} \in A$. Using $M_d \su A$ again we find
$x^d y^s \in A$, which is impossible.

Next suppose $e>d$; then $m\ge e$. Turning to the polynomial
$\tfrac{\del ^{e-d}}{\del x^{e-d}}F$ we reduce this case to the case when
$e=d$. \hfill $\square$

If $d=0$, then it follows from Lemma \ref{l4} that if $F\in A$ and $[F]_0 =
\ldots = [F]_{s-1} =0$, then $F=0$. That is, if $d=0$ then $A$ is an L-module of
order $s$. Then $A=M_0 +A$ gives a representation needed.
Therefore, we may assume that $d\ge 1$.

If $s=0$, then it follows from Lemma \ref{l4} that if $F\in A$, then
$\deg [F]_n <d$ for every $n$. Thus $A\su M_d$ and, consequently, we have
$A=M_d$. Putting $M =\{ 0\}$ (which is an $L$-module of arbitrary order with an
arbitrary $L$), we obtain $A= M_d +M$. Therefore, we may assume $s\ge 1$.

\begin{notation}
{\rm If $\phi \in \cc [x]^s$ and $\phi =(f_0 \stb f_{s-1} )$, then we use the
notation $\phi ' =(f'_0 \stb f'_{s-1} )$. We say that a subset $V$ of $\cc [x]^s$
is closed under differentiation, if $\phi \in V$ implies $\phi ' \in V$.
Note that $V_{A,s}$ is closed under differentiation by \eqref{e16}.

Let $V_k =\{ (f_1 \stb f_s )\in V_{A,s} \colon \deg f_i <k \ (i=1\stb s)\}$
for every integer $k$. Note that $V_k$ is also closed under differentiation.
We have $V_d = \{ (f_1 \stb f_s )\in \cc [x]^s \colon \deg f_i <d \
(i=1\stb s)\}$, as $M_d \su A$.

For every polynomial $F\in \cc [x,y]$ we denote
$$\Phi (F)=( [F]_0 \stb [F]_{s-1} )\in \cc [x]^s .$$
Clearly, $\Phi$ is a linear map from $\cc [x,y]$ onto $\cc [x]^s$, and maps
$A$ onto $V_{A,s}$. 

It follows that there exists a linear map $\phi \mapsto F_\phi$ from
$V_{A,s}$ into $A$ such that $\Phi (F_\phi )=\phi$ for every
$\phi \in V_{A,s}$. 
}
\end{notation}
\begin{lemma}\label{l5}
For every integer $k>d$ there is a linear map $L \colon V_k \to \cc [x]$
with the following properties.
\begin{enumerate}[{\rm (i)}]
\item For every $\phi \in V_k$ we have
\begin{equation}\label{e30}
\deg (L(\phi )- [F_\phi ]_s )<d.
\end{equation}
\item $L(\phi ')=L(\phi )'$ for every $\phi \in V_k$.
\item $L(\phi )=0$ for every $\phi \in V_d$.
\end{enumerate}
\end{lemma}

\noi
{\it Proof.} Let $X$ denote the quotient space of the linear space $V_k$
modulo the linear subspace $V_d$. (That is, let $X=V_k /V_d$.)
Since the linear space $V_k$ is of finite dimension (its dimension is at
most $k^s$), so is $X$. Let $\phi \mapsto \ol \phi$ denote the natural
homomorphism from $V_k$ into $X$. That is, let $\ol \phi = \phi +V_d$
for every $\phi \in V_k$.

The derivation $\phi \mapsto \phi '$ maps $V_d$ into itself. Therefore, we
can define the derivation on $X$ by $D(\ol \phi )=\ol{\phi '}$
$(\phi \in V_k )$. 

It is clear that $D$ is a nilpotent linear map from $X$ into
itself. By \cite[\S 57, Theorem 2, p.~111]{H}, there are positive integers $r,
q_1 \stb q_r$ and elements $u_1 \stb u_r \in X$ such that
$D^{q_i} u_i =0$ for every $i=1\stb r$, and the elements $D^{j} u_i$
$(i=1\stb r,\ j=0\stb q_i -1)$ form a basis for $X$.
Let $\psi _1 \stb \psi _r \in V_k$ be such that $u_i =\ol{\psi _i}$
$(i=1\stb r)$. We put
\begin{equation}\label{e25}
\La (D^j u_i ) =([F_{\psi_i}]_s )^{(j)}
\end{equation}
for every $i=1\stb r$ and $j=0\stb q_i -1$, and extend $\La$ linearly to $X$.
We define $L (\phi )=\La (\ol \phi )$ for
every $\phi \in V_k$. Then $L\colon V_k \to \cc [x]$ is linear. We show that $L$
has properties (i)-(iii).

If $\phi \in V_d$, then $\ol \phi =0$, $L(\phi )=\La (\ol \phi )=0$,
and thus (iii) holds.

Next we prove (i). Since $L$ and the map $\phi \mapsto F_\phi$ are both
linear, the set of elements $\phi \in V_k$ satisfying \eqref{e30} is a
linear subspace of $V_k$. 
Therefore, in order to prove (i) it is enough to check that \eqref{e30} holds
for a set of polynomials generating $V_k$. We show that $\Psi =\{ \psi _i^{(j)}
\colon i=1\stb r, \ j=1\stb q_i -1 \} \cup V_d$ is such a set. Indeed,
let $\phi \in V_k$. Since $D^{j} u_i$ is a basis for $X$, we have
\begin{equation}\label{e26}
\ol \phi = \sumir \sum_{j=1}^{q_j} \la _{i,j} D^j u_i
\end{equation}
with suitable complex coefficients $\la _{i,j}$.
Now $\ol{\al ' } =D \ol \al$ $(\al \in V_k )$ implies that
$\ol{\psi _i^{(j)}} =D^j u_i$ for every $i,j$, and thus the right hand side of
\eqref{e26} equals the image under the natural homomorphism of a linear
combination of the elements $\psi _i^{(j)}$. Thus the difference of $\phi$ and
this linear combination belongs to $V_d$, showing that $\Psi$ generates $V_k$.

If $\phi \in V_d$, then $L(\phi )=0$ by (iii) and $\deg [F_\ph ]_s <d$
by Lemma \ref{l4}, and thus \eqref{e30} holds.

If $\phi =\psi _i ^{(j)}$, then we have $\Phi (F_{\psi _i})=\psi _i$,
$$\Phi (\tfrac{\del ^j}{\del x^j} F_{\psi _i} )= \psi _i ^{(j)} =
\Phi ( F_{\psi _i ^{(j)}}) =\Phi (F_\phi ),$$
and thus $\Phi (\tfrac{\del ^j}{\del x^j} F_{\psi _i} - F_{\phi}) =0$. In other
words, the first $s$ coordinate polynomials of $\tfrac{\del ^j}{\del x^j}
F_{\psi _i} - F_{\phi}$ are zero. By Lemma \ref{l4} it follows that
\begin{align*}
  d &>\deg ([\tfrac{\del ^j}{\del x^j} F_{\psi _i}  - F_{\phi} ]_s)=
\deg ([\tfrac{\del ^j}{\del x^j} F_{\psi _i} ]_s - [F_{\phi}]_s)\\
&=  \deg (([F_{\psi_i}]_s )^{(j)}  - [F_{\phi} ]_s) =
  \deg (\La (D^j u_i ) -[F_\phi ]_s ) \\
&  =\deg (L(\phi )-[F_\phi ]_s ),
\end{align*}
which proves (i).

We turn to the proof of (ii). Since $L$ is linear and $L(\phi )=0$ if
$\phi \in V_d$, in order to prove (ii) it is enough to check that
$L(\phi ')=L(\phi )'$  holds in the cases when $\phi =\psi _i ^{(j)}$.
Let $1\le i\le r$ and $0\le j\le q_i -1$ be fixed. If $j<q_i -1$,
then $L(\phi )' = [F_{\psi_i}]_s ^{(j+1)} =L(\phi ')$, and we are done.

If $j=q_i -1$, then $\phi ' =\psi _i ^{(q_i )} =0$, so we have
$\phi \in V_1 \su V_d$. Then $L(\phi )=0$ by (iii), and $L(\phi ')=
L(\phi )'=0$ follows. \hfill $\square$

\begin{lemma}\label{l6}
There exists a linear map $L\colon V_{A,s} \to \cc [x]$ such that
\begin{enumerate}[{\rm (i)}]
\item \eqref{e30} holds for every $\phi \in V_{A,s}$,
\item $L(\phi ')=L(\phi )'$ for every $\phi \in V_{A,s}$, and
\item $L(\phi )=0$ for every $\phi \in V_d$.
\end{enumerate}
\end{lemma}

\noi
{\it Proof.} 
We define
\begin{equation}\label{e28}
L(\phi )=[F_{\phi}]_s + u_{\phi ,d-1} x^{d-1} +\ldots +u_{\phi ,1} x +u_{\phi ,0} ,
\end{equation}
where $u_{\phi ,i}$ is an unknown for every $\phi \in V_{A,s}$ and $i=0\stb d-1$.
We show that we can assign values to these unknowns in such a way that the
resulting map $L$ satisfies the requirements. Since the map $\phi
\mapsto F_\phi$ is linear, $L$ will be linear if
\begin{align*}
 u_{\la \phi +\mu \psi ,d-1} & x^{d-1} +\ldots +
u_{\la \phi +\mu \psi ,1} x +  u_{\la \phi +\mu \psi ,0}=\\
&\la (u_{ \phi ,d-1} x^{d-1} +\ldots +u_{ \phi ,1} x
  +u_{ \phi ,0}) + \\
&\mu (u_{ \psi ,d-1} x^{d-1} +\ldots +u_{ \psi ,1} x
+u_{ \psi ,0})
\end{align*}
holds for every $\phi ,\psi \in V_{A,s}$ and $\la ,\mu \in \cc$. 
It is clear that condition (i) is satisfied with any choice of the unknowns
$u_{ \phi ,i}$. Condition (ii) is satisfied if
\begin{align*}
&  [F_{\phi '}]_s +  u_{\phi ' ,d-1} x^{d-1} +\ldots + u_{\phi ' ,1} x +  u_{\phi ' ,0}=\\
  & ([F_{\phi}]_s )' + (d-1) u_{\phi ,d-1} x^{d-2} +\ldots +u_{ \phi ,1} 
\end{align*}
holds for every $\phi \in V_{A,s}$. Finally, (iii) is satisfied if the right
hand side of \eqref{e28} is zero for every $\phi \in V_d$. 
Summing up: in order that $L$ satisfy the conditions, the unknowns
$u_{ \phi ,i}$ must satisfy a certain infinite system of linear equations $S$.
We have to show that $S$ is solvable.
It is well-known that a system $S$ of linear equations is solvable if and only
if every finite subsystem of $S$ is solvable. Now a finite subsystem $T$ of $S$
only involves a finite number of elements $\phi \in V_{A,s}$. Then there is a $k$
such that all these elements belong to $V_k$. As we proved above, there is a
map $L$ on $V_k$ satisfying (i)-(iii) on $V_k$.

Now condition (i) implies that $L$ is of the form \eqref{e28}
with concrete values of the unknowns $u_{ \phi ,i}$ for every $\phi \in V_k$.
These values constitute a solution of the subsystem $T$, showing that $T$
is solvable. Therefore, $S$ is solvable, proving the existence of $L$ with
the required properties. \hfill $\square$

\noi
{\it Proof of Theorem \ref{t2}.}
Fix a map $L$ as in Lemma \ref{l6}. We prove that if $\phi =
(f_{0} \stb f_{s-1}) \in V_{A,s}$, then the recursion $f_n = L(f_{n-s} \stb f_{n-1})$
$(n=s,s+1,\ldots )$ defines a sequence of polynomials such that $f_n =0$
for every $n$ large enough, and
\begin{equation}\label{e31}
\deg (f_n -[F_\phi ]_n )<d
\end{equation}
for every $n$. It is clear that \eqref{e31} holds for every $n<s$.

Let $k\ge s$, and suppose we have defined $f_n$ for every $n<k$ such that
\eqref{e31} holds for every $n<k$. Let $\psi =(f_{k-s}\stb f_{k-1})$. We have
$F_\phi \in A$ and $G=\tfrac{\del ^{k-s}}{\del y^{k-s}} F_\phi \in A$.
Since $\deg (f_{k-s+i} -[F_\phi ]_{k-s+i} )<d$ and $[G]_i =[F_\phi ]_{k-s+i}$ for
every $i<s$, we have $\Phi (G)-\psi \in V_d$. Since $V_d \su V_{A,s}$ and
$\Phi (G)\in V_{A,s}$, we obtain $\psi \in V_{A,s}$. Therefore,
$L(f_{k-s}\stb f_{k-1} )$ is defined. Let $f_k =L(f_{k-s}\stb f_{k-1})$. By (i) of
Lemma \ref{l6}, we have $\deg (f_k - [F_\psi ]_s )<d$. 

Now $G-F_\psi \in A$ and $\Phi (G-F_\psi ) \in V_d$. By Lemma \ref{l4},
this implies $\deg [G-F_\psi ]_s <d$. Since $[G]_s =[F_\phi ]_k$, we obtain
$\deg ([F_\phi ]_k -[F_\psi ]_s )<d$ and
$\deg (f_k - [F_\phi ]_k )<d$. This proves that the recursion
$f_n = L(f_{n-s} \stb f_{n-1})$ defines $f_n$ for every $n$ such that
\eqref{e31} holds for every $n$.

Since $F_\phi \in A$, there is an $N$ such that $[F_\phi ]_n =0$ for every $n\ge
N$. Then $\deg f_n <d$ for every $n\ge N$. If $n>N+s$, then
$(f_{n-s} \stb f_{n-1}) \in V_d$ by \eqref{e31}, and thus 
$f_n =L(f_{n-s} \stb f_{n-1})=0$ by (iii) of Lemma \ref{l6}.
Therefore $f_n =0$ for every $n$ large enough. Let $H_\phi$ denote the polynomial
$\sum_{n=0}^\infty f_n \cd \frac{y^n}{n!}$. Then \eqref{e31} implies that
$H_\phi -F_\phi \in M_d$. Since $M_d \su A$ and $F_\phi \in A$, it follows
that $H_\phi \in A$.

Let $M$ be the set of polynomials $H_\phi$, where $\phi \in V_{A,s}$.
Then we have $M\su A$. It is easy to see that the map $\phi \mapsto H_\phi$
is linear, and thus $M$ is a linear subspace of $A$.
It is also easy to check that $F\in  M$ implies $\tfrac{\del}{\del y}F\in M$.
Now (ii) of Lemma \ref{l6} implies that $\tfrac{\del}{\del x}F_\phi =H_{\phi '}
\in M$ for every $\phi \in V_{A,s}$. Thus $M$ is also closed under partial
differentiation w.r.t. $x$. Consequently, $M$ is a module. It is clear
that $M$ is an L-module of order $s$.

If $F\in A$, then $\phi =\Phi (F)\in V_{A,s}$. Now $F-F_\phi \in M_d$
by Lemma \ref{l4}, and $H_\phi -F_\phi \in M_d$ by \eqref{e31}. Thus
$$F=((F-F_\phi )-(H_\phi -F_\phi )) +H_\phi \in M_d +M,$$
which proves $A=M_d +M$. \hfill $\square$

\begin{remark} \label{r1}
{\rm In the representation $A=M_d +M$, where $M$ is an L-module, the value of
$d$ is unique (see Corollary \ref{c1}). However, the term $M$ is not unique in 
general, as the following example shows.

Let $M_{d,s}$ denote the set of polynomials of the form \eqref{e1}, where
$\deg f_n <d$ for every $n<s$, and $f_n =0$ for every $n\ge s$. It is clear
that $M_{d,s}$ is an L-module of order $s$. Since $M_{d,s} \su M_d$, we have
$M_d =M_d +\{ 0\} =M_d +M_{d,s}$.

We can make the representations unique if we restrict the L-module terms.
Note that the proof of Theorem \ref{t2} produces 
L-modules with a linear map $L$ such that
$L(f_1 \stb f_s )=0$ whenever $\deg f_i <d$ $(i=1\stb s)$; see (iii) of Lemma
\ref{l6}.
We may also assume that $M_{d,s} \su M$, since otherwise we replace $M$ by
$M+M_{d,s}$. Now, it is easy to check that the representation $A=M_d +M$
is unique, if we require that the L-module $M$ should satisfy both $M_{d,s} \su
M$ and $L(f_1 \stb f_s )=0$ whenever $\deg f_i <d$ $(i=1\stb s)$.}
\end{remark}

\section{Indecomposable submodules} \label{s1}

\begin{definition} 
{\rm We say that a submodule $M$ of $\cc [x_1 \stb x_n ]$ is {\it decomposable},
if $M$ can be represented as the sum of finitely many proper submodules of $M$.
Otherwise the submodule $M$ is {\it indecomposable}.
}
\end{definition}

\begin{proposition}\label{p1}
Every submodule of $\cc [x_1 \stb x_n ]$ is the sum of finitely many
indecomposable submodules.
\end{proposition}

\noi
{\it Proof.} The family of submodules of $\cc [x_1 \stb x_n ]$ has the minimal
condition; that is, if $M_1 \sp M_2 \sp \ldots$ are submodules of
$\cc [x_1 \stb x_n ]$, then there is a positive integer $K$ such that $M_k =M_K$
for every $k \ge K$ (see \cite[Lemma 8]{LSz}). Therefore, if the statement of
the proposition is not true, then there is a minimal counterexample $M$. Then
$M$ must be decomposable. If $M=A_1 +\ldots +A_k$, where $A_1 \stb A_k$ are 
proper submodules of $M$ then, by the minimality of $M$, each $A_i$ is the sum
of finitely many indecomposable submodules. Then the same is true for $M$, which
is impossible. \hfill $\square$

It is not clear if the representation of a module as the sum of indecomposable
submodules containing a minimal number of terms is unique or not.

In the following we confine ourselves to the submodules of $\cc [x,y]$ (except
in Remark \ref{r2}).
It follows from Theorem \ref{t2} that {\it if $M$ is an indecomposable submodule
of $\cc [x,y]$, then either $M=M_d$ for some $d$ or $M$ is an L-module of order
$s$ for some $s$.}

Our next aim is to show that $M_d$ is indecomposable for every $d$, and so is
every L-module of order $1$.
\begin{lemma}\label{l3}
The system of translation invariant linear subspaces of $\cc [x]^s$ has the
minimal condition.
\end{lemma}

\noi
{\it Proof.} We prove the statement by induction on $s$. Since every
translation invariant linear subspace of $\cc [x]$ equals $\cc [x]$ or
$\{ f\in \cc [x]\colon \deg f <d\}$ for some $d\ge 0$, it easily follows that
the statement is true for $s=1$.

Let $s\ge 1$, suppose that the statement is true for $s$, and let
$V_1 \sp V_2 \sp \ldots$ be translation invariant linear subspaces of $\cc [x]^{s+1}$. We have to show that $V_n =V_{n+1}=\ldots$ if $n$ is large enough.

Put $A_n =\{ f\in \cc [x]\colon (0\stb 0,f)\in V_n \}$ for every $n$.
Since $A_n$ is a translation invariant linear subspace of $\cc [x]$ and
$A_1 \sp A_2 \sp \ldots$, there is an $N_1$ such that $A_n =A_{N_1}$ for every
$n\ge N_1$. Let 
$$B_n =\{ (f_1 \stb f_s )\in \cc [x]^s \colon \exists \ f, \ 
(f_1 \stb f_s ,f)\in V_n \} .$$
Then $B_n$ is a translation invariant linear subspace of $\cc [x]^s$ and
$B_1 \sp B_2 \sp \ldots$. By the induction hypothesis it follows that 
there is an $N_2$ such that $B_n =B_{N_2}$ for every $n\ge N_2$. Let 
$N=\max (N_1 ,N_2 )$; we prove that $V_n =V_{N}$ for every $n\ge N$.
Let $n\ge N$ and $(f_1 \stb f_{s+1} )\in V_N$ be given; we prove
$(f_1 \stb f_{s+1} )\in V_n$.

We have $(f_1 \stb f_s )\in B_N =B_n$, and thus there is a $g$ such that
$(f_1 \stb f_s ,g)\in V_n \su V_N$. From $(f_1 \stb f_{s+1} )\in V_N$ we obtain
$(0\stb 0, f_{s+1} -g)\in V_N$, $f_{s+1} -g\in A_N =A_n$ and $(0\stb 0, f_{s+1} -g)
\in V_n$ . Thus
$$(f_1 \stb f_{s+1} )=(f_1 \stb f_s ,g) + (0\stb 0, f_{s+1} -g) \in V_n ,$$
and the proof is complete. \hfill $\square$

\begin{theorem} \label{t3}
If $A,B \su \cc [x,y]$ are L-modules, then so is $A+B$.
\end{theorem}

\noi
{\it Proof.} Suppose $A$ is of order $s_1$ and $B$ is of order $s_2$. If
$s=\max(s_1 ,s_2 )$, then $A,B$ are both of order $s$. For every $k>s$
we denote by $Z_k$ the set of $s$-tuples $(f_1 \stb f_s )\in V_{A,s} \cap V_{B,s}$
such that if $F\in A$, $G\in B$ and $[F]_n =[G]_n =f_n$ for every $n<s$, then 
$[F]_n =[G]_n$ for every $n< k$.

It is easy to check that $Z_k$ is a translation invariant linear subspace of
$\cc [x]^s$, and $Z_{s+1} \sp Z_{s+2} \sp \ldots$. By Lemma \ref{l3}, there is a
$K > s$ such that $Z_k =Z_K$ for every $k\ge K$. We prove that $A+B$ is an
L-module of order $K$.

Let $S \in A+B$ such that $[S]_n =0$ for every $n<K$. We show that $[S]_K =0$.
Let $S=F+G$, where $F\in A$ and $G\in B$. Then $[F]_n +[G]_n =0$; that is,
$[F]_n =-[G]_n$ for every $n<K$. Since $-G\in B$, it follows that $([F]_0 \stb
[F]_{s-1}) \in V_{A,s} \cap V_{B,s}$. We prove $([F]_0 \stb [F]_{s-1}) \in Z_K$.

Suppose $H\in A$ and $[H]_n =[F]_n$ for every $n<s$.
Since $F,H\in A$ and $A$ is an L-module of order $s$, it follows that $F=H$, and
thus $[H]_n =[F]_n$ for every $n<K$. Similarly, if $P\in B$ and
$[P]_n =[F]_n =-[G]_n$ for every $n<s$, then $P=-G$, and
thus $[P]_n =-[G]_n =[F]_n$ for every $n<K$, proving
$([F]_0 \stb [F]_{s-1}) \in Z_K$.

Since $Z_K =Z_{K+1}$, we find $([F]_0 \stb [F]_{s-1}) \in Z_{K+1}$. This implies
$[F]_K =-[G]_K$; that is, $[S]_K=0$. \hfill $\square$

\begin{remark}
{\rm The proof above does not give any estimate of the order of $A+B$.
We do not know if the order of $A+B$ is bounded from above by, say, the sum
of the order of $A$ and of $B$.}
\end{remark}

\begin{theorem} \label{t4}
$M_d$ is indecomposable for every $d$.
\end{theorem}

\noi
{\it Proof.} Suppose this is not true, and let $M_d =A_1 +\ldots +A_k$,
where $A_1 \stb A_k$ are proper submodules of $M_d$. By Theorem \ref{t2}, we
have $A_i =M_{d_i} +B_i$, where $B_i$ is an L-module for every $i$. By Theorem
\ref{t3} we find that $B=B_1 +\ldots +B_k$ is an L-module. It is clear that
$M_{d_1} +\ldots +M_{d_k} =M_e$, where $e=\max_{1\le i\le k} d_i$. Therefore,
$M_d =A_1 +\ldots +A_k$ gives $M_d =M_e +B$, where $B$ is an L-module. By
Corollary \ref{c1}, we have $e=d$, and thus
$d_i =d$ for a suitable $1\le i\le k$. Then $M_d = M_{d_i} \su A_i$, which
is impossible, since $A_i$ is a proper submodule of $M_d$. \hfill $\square$

\begin{theorem} \label{t5}
Every L-module of order $1$ is indecomposable.
\end{theorem}

\noi
{\it Proof.} Let $M$ be an L-module of order $1$, and suppose $M =A_1 +\ldots
+A_k$, where $A_1 \stb A_k$ are proper submodules of $M$. Each of the linear
spaces $V_{M,1}$ and $V_{A_i ,1}$ $(i=1\stb k)$ equals one of $\cc [x]$ or
$\{ f\in \cc [x]\colon \deg f <d \}$ for some $d\ge 0$. Since
$$V_{M,1} = V_{A_i ,1} +\ldots +V_{A_k ,1} ,$$
it follows that $V_{M,1} = V_{A_i ,1}$ for a suitable $i$. Then we have $M=A_i$
by the definition of L-modules, which is impossible. \hfill $\square$

\begin{remarks}
{\rm (i) There are decomposable L-modules: if $A, B$ are L-modules,
$A\subsetneq B$ and $B\subsetneq A$, then $A+B$ is a decomposable L-module.

\noi
(ii) There are indecomposable L-modules of order $>1$. Indeed, let
$$A^* =\{ f(y,x)\colon f(x,y)\in A\}$$
for every $A\su \cc [x,y]$. It is clear that if $M$ is a module, then so is
$M^*$, and if $M$ is indecomposable then so is $M^*$. Thus $M_2^*$ is
indecomposable. On the other hand, $M_2^* =\{ f(x)+g(x)y\colon f,g \in
\cc [x] \}$. It is clear that $M_2^*$ is an L-module of order $2$.
Similarly, $M_d^*$ is an indecomposable L-module of order $d$ for every $d$.

\noi
(iii) It follows from the definition that a submodule of an L-module is also
an L-module. Also, we have $M_d \not\subset M$ for every $d\ge 1$ and for every
L-module $M$. Indeed, if $M$ is an L-module of order $s$, then $y^s \notin M$
and $y^s \in M_1 \su M_d$. From these observations it follows that}
the representation of the submodules of $\cc [x,y]$ as sums of indecomposable
submodules containing a minimal number of terms is unique if and only if this is
true for L-modules.
\end{remarks}

\begin{remark} \label{r2}
{\rm
We show that Theorem \ref{t2} does not have a straightforward generalization to
$\cc [x,y,z]$. Such a generalization would operate with modules defined as
follows. Generalizations of the modules $M_d$ could be defined as the set of
polynomials
\begin{equation}\label{e17}
F(x,y,z)=\sum_{n=0}^\infty f_n (x,y)\frac{z^n}{n!}
\end{equation}
such that the degrees of the polynomials $f_n (x,y) \in \cc [x,y]$ satisfy some 
prescribed inequalities. Let these modules be called bounded modules.
We call $M$ an L-module of order $s$ if, whenever the polynomial in \eqref{e17}
belongs to $M$ and $f_n =0$ for every $n<s$, then $F=0$.

Now suppose Theorem \ref{t2} had a generalization to $\cc [x,y,z]$. It would
claim that every proper submodule of $\cc [x,y,z]$ is the sum
of bounded modules and L-modules. Then it would follow that whenever $M$ is an
indecomposable submodule of $\cc [x,y,z]$, then either $M$ is bounded, or $M$ is
an L-module.

We show that this is false, no matter how we define bounded modules. Let $M$
be the set of polynomials
\begin{equation}\label{e18}
\sum_{n=0}^\infty (a_n (x+y)+b_n ) \cd \frac{z^n}{n!} ,
\end{equation}
where $a_n$ and $b_n$ are complex numbers and $a_n =b_n =0$ if $n$ is large
enough. It is clear that $M$ is a
module. Suppose $M=A_1 +\ldots +A_k$, where $A_1 \stb A_k$ are submodules of
$M$. Then there is an $i$ such that $A_i$ has the following property: for
every $N$ there is a polynomial of the form \eqref{e18} belonging to $A_i$ and
such that $a_n \ne 0$ for at least one index $n>N$. It is easy to check
that this condition implies $A_i =M$, and thus $M$ is indecomposable.

Now $M$ is not a bounded module, since no matter how we prescribe the
inequalities satisfied by the elements of $M$, the polynomial $x+2y$ would
also satisfy these conditions, but $x+2y\notin M$. It is also clear that
$M$ is not an L-module, since the coordinate polynomials $a_n (x+y)+b_n$ can be
chosen independently. This shows that no generalization of the form described
above is possible.}
\end{remark}

\section{Construction of a submodule of $\cc [x,y]$ which is not closed}

We equip $\cc [x_1 \stb x_n ]$ with the topology of uniform convergence on
compact sets. The closure of a set $M\su \cc [x_1 \stb x_n ]$ w.r.t. this
topology is denoted by $\cl M$.

\begin{theorem}\label{t6}
There exists a module $M\su \cc [x,y]$ such that $x\in \cl M$ but $x\notin M$.
\end{theorem}

\noi
{\it Proof.}
We use the notation $e(x)=e^x$, $e_2 (x)=e(e(x))$ and $e_3 (x)=e(e_2 (x))$.
Then we have 
\begin{equation}\label{e14}
e_3 (n-1)^{e(n)}/e_3 (n)\to 0 \ (n\to \infty ).
\end{equation}
Indeed, we have $n+e(n-1)-e(n)\to -\infty$, hence $e(n)\cd e_2 (n-1)/e_2 (n)
\to 0$, hence $e(n)\cd e_2 (n-1) -e_2 (n) \to -\infty$, hence
$e_3 (n-1)^{e(n)}/e_3 (n)\to 0$ as $n\to \infty$.

Let $a_1 =1$ and $a_n =-e_3 (n)$ for every
$n\ge 2$. We put $L(f)=a_1 f' +a_2 f'' +\ldots$ for every $f\in \cc [x]$.
Note that the number of nonzero terms in the sum is finite for every $f\in
\cc [x]$. Let $M$ denote the set of polynomials
\begin{equation}\label{e9}
\sum_{n=0}^\infty L^n (f) (x)\cd \frac{y^n}{n!},
\end{equation}
where $f\in \cc [x]$ is arbitrary. By Lemma \ref{l1}, $M$ is a submodule of
$\cc [x,y]$. Now we have $x\notin M$. Indeed, if $F(x,y)$ is defined by
\eqref{e1} and $F(x,y)=x$, then $f_0 (x)=x$ and $f_1 =0$. However, we have
$f_1 = L (f_0 ) = L (x)=1$, a contradiction. (In fact, the same argument
gives $f\notin M$ for every $f\in \cc [x]$ with $\deg f \ge 1$.)

We put $\ep _n =|a_n \cd n!| \am$, $g_n (x)=x +\ep _n x^n$ and
$G_n (x,y)=\sum_{k=0}^n L^k (g_n )\cd \tfrac{y^k}{k!}$ for every $n$.
We show that the sequence of polynomials $G_n$ converges to $x$ locally
uniformly on $\cc ^2$. Since $G_n \in M$ for every $n$, this will prove that 
$x\in \cl M$.

If $f\in \cc [x]$, $f=\sum_{i=0}^n c_i x^i$, then we put $\V f\V =
\max_{0\le i\le n} |c_i |$. Clearly, $\V .\V$ is a norm on $\cc [x]$.
If $\deg f \le n$ and $|x|\le e^n$, then
\begin{equation}\label{e11}
\begin{split}
|f(x)|&\le \V f\V \cd (1+e (n) +e(2n)+\ldots +e (n^2 ))<
(n+1)\cd e (n^2 )\cd \V f\V\\
& < e_2 (n)\cd \V f\V <e_3 (n-1)\cd \V f\V
\end{split}
\end{equation}
if $n>n_0$. If $f\in \cc [x]$ and $\deg f \le n$, then $\V f'\V \le n\cd
\V f\V$. Therefore, 
\begin{equation}\label{e10}
\begin{split}
  \V L f\V &\le \sum_{i=1}^n |a_i | \cd \V f^{(i)} \V \le n! \sum_{i=1}^n |a_i | \cd \V f\V \\
&\le n! \cd n \cd e_3 (n)  \cd \V f\V \\
&\le e_2 (n)\cd e_3 (n)< e_3 (n)^2 \cd \V f\V
\end{split}
\end{equation}
for every $f\in \cc [x]$ with $\deg f \le n$. (We used the trivial estimate
$n! \cd n\le n^n <e(n^2 )<e_2 (n)$.) 
Let $n\ge 2$ be fixed. If $|x|\le e^n$, then
\begin{equation}\label{e12}
|g_n (x)-x| \le \ep _n \cd e (n^2) < e(n^2)/ e_3 (n) =e(n^2 -e_2 (n) ).
\end{equation}
Now we have
\begin{align*}
L (g_n ) &= g_n' +\sum_{i=2}^{n-1} a_i g_n^{(i)} +a_n g_n^{(n)}\\
&  =1+n\cd \ep _n  x^{n-1} + \sum_{i=2}^{n-1} a_i \cd n(n-1)\cdots (n-i+1) \ep _n
x^{n-i} +a_n \cd n! \cd \ep _n \\
&= n\cd \ep _n  x^{n-1} + \sum_{i=2}^{n-1} a_i \cd n(n-1)\cdots (n-i+1) \ep _n
x^{n-i},
\end{align*}
and thus
$$\V L (g_n ) \V \le \ep _n \cd \max_{1\le i \le n-1} |a_i | \cd n! < \ep _n \cd
e_3 (n-1)\cd n! =e_3 (n-1)/e_3 (n).$$
Since $\deg L (g_n ) =n-1$, \eqref{e10} gives
$$
\V L ^k (g_n ) \V \le e_3 (n-1)^{2k-2} \cd \V L g_n \V < e_3 (n-1)^{2k-1} /e_3 (n).
$$
for every $k\ge 2$. Then we find, by \eqref{e11}, that if $2\le k\le n$ and
$|x|\le e^n$, then
\begin{equation}\label{e13}
  |L^k (g_n ) (x)|\le e_3 (n-1)^{2n}/e_3 (n) .
\end{equation}
If $|x|\le e^n$ and $|y|\le e^n$, then it follows from \eqref{e12} and
\eqref{e13} that
$$|G_n (x,y)-x|<e(n^2 -e_2 (n) ) +n\cd e(n^2)\cd e_3 (n-1)^{2n}/e_3 (n)$$
if $n>n_0$. Since $e(n^2 -e_2 (n) ) \to 0$ and
$$e(n^2)\cd e_3 (n-1)^{2n}<e_3 (n-1)^{2n+1} <e_3 (n-1)^{e(n)},$$
it follows from \eqref{e14} that $G_n (x,y)\to x$ locally uniformly on
$\cc ^2$. \hfill $\square$

\begin{remarks}
{\rm (i) It is easy to see that $\cc [x_1 ,x_2 ]$ is a closed submodule of
$\cc [x_1 \stb x_n ]$ for every $n\ge 2$. Therefore, Theorem \ref{t6} implies
that} for every $n\ge 2$ there exists a submodule of $\cc [x_1 \stb x_n ]$
which is not closed.

\noi
{\rm (ii)  
Using an elaborate version of the proof of Theorem \ref{t6} one can show that}
{\it there are L-modules of order $1$ which are everywhere dense in
$\cc [x,y]$ w.r.t. the topology of uniform convergence on compact sets.}
\end{remarks}

\noi
R\'enyi Institute, Budapest, Hungary

\noi
ELTE E\"otv\"os Lor\'and University, Budapest, Hungary

\noi
e-mail: {\tt kigergo57@gmail.com, miklos.laczkovich@gmail.com}

\end{document}